\documentclass{amsart}

\renewcommand{\epsilon}{\varepsilon}
\newcommand{\varkappa}{\kappa}

\newcommand{\carea}{c_{\mathit{area}}}

\newcommand{\Z}{{\mathbb Z}}
\newcommand{\Q}{{\mathbb Q}}
\newcommand{\R}{{\mathbb R}}
\newcommand{\C}{{\mathbb C}}

\newcommand{\Vol}{\operatorname{Vol}}

\newcommand{\SL}{\operatorname{SL}(2,\R)}
\newcommand{\GL}{\operatorname{GL}^+(2,\R)}

\newlength{\halfbls}
\newcommand{\cH}{{\mathcal H}}
\newcommand{\cM}{{\mathcal M}}

\newtheorem{MainConjecture}{Main Conjecture}

\newtheorem{ConditionalCorollary}{Conditional Corollary}
\newtheorem{Conjecture}{Conjecture}

\theoremstyle{remark}
\newtheorem{Remark}{Remark}

\begin{document}
\begin{picture}(0,0)(0,0)
\put(-20,60){to appear in Arnold Mathematical Journal}
\put(-20,50){DOI 10.1007/s40598-015-0028-0}
\end{picture}

\author{Alex Eskin}
\thanks{Research  of  the first author is partially supported  by
NSF grant.}
\address{
Department of Mathematics,
University of Chicago,
Chicago, Illinois 60637, USA\\
}
\email{eskin@math.uchicago.edu}

\author{Anton Zorich}
\thanks{Research  of  the second author is partially supported  by
IUF and by ANR}
\address{
Institut de math\'ematiques de Jussieu,
Institut Universitaire de France,
Universit\'e Paris 7, France}
\email{Anton.Zorich@gmail.com}

\dedicatory{To the memory of V.~I.~Arnold}

\date{August 31, 2015}

\title[Volumes and Siegel--Veech constants in large genera]
{Volumes of strata of Abelian differentials and
Siegel--Veech constants in large genera}

\begin{abstract}
We  state  conjectures  on  the asymptotic behavior of the volumes of
moduli  spaces  of  Abelian  differentials  and  their  Siegel--Veech
constants  as  genus  tends to infinity. We provide certain numerical
evidence,  describe  recent advances and the state of the art towards
proving these conjectures.
\end{abstract}

\maketitle

\section*{Introduction}

The  first  efficient approach for computing volumes of the strata in
the  moduli spaces of Abelian differentials was found by A.~Eskin and
A.~Okounkov~\cite{Eskin:Okounkov}   about   fifteen  years  ago.  The
algorithm  was implemented by A.~Eskin in a rather efficient computer
code  which  already  at  this time allowed to compute volumes of all
strata up to genus 10, and volumes of some strata, like the principal
one, up to genus 60 (or more).

About   the  same  time  together  with  H.~Masur  we  expressed  the
Siegel--Veech   constant   $\carea$   of   any   stratum  of  Abelian
differentials  in  terms  of an explicit polynomial in volumes of
the  ``principal  boundary''  strata devided by the volume of the the
original  stratum,  see~\cite{Eskin:Masur:Zorich}.  The  sum  of  the
Lyapunov  exponents  of  the  Hodge  bundle  over  the  Teichm\"uller
geodesic  flow  can  be expressed by a simple formula in terms of the
Siegel--Veech  constant $\carea$ (this formula, is recently proved by
M.~Kontsevich  and  the  authors  in~\cite{Eskin:Kontsevich:Zorich}).

Direct  computations of volumes of the strata in genera accessible at
this  time  combined  with  numerical  experiments  with the Lyapunov
exponents  providing  an  approximate  value of $\carea$ inspired two
conjectures  on asymptotic behavior of volumes of the strata and of
the  corresponding Siegel--Veech constant $\carea$; these conjectures
were stated at the end of 2003.

During the last decade we have shared this conjecture with colleagues
who  were  interested in it, but the fact that it was never published
lead  to  some misinterpretations due to misundersting in conventions
on normalization. Now, when the formula\cite{Eskin:Kontsevich:Zorich}
for  the  sum  of  the  Lyapunov  exponents  is  finally  proved, the
conjectures  become  much  more  important. Besides, there are recent
encouraging  advances with these conjectures, so it is an occasion to
state the conjectures and describe the state of the art.

\section{Background material and motivations}

Holomorphic  $1$-forms  on  a  given Riemann surface $C$ of genus $g$
form  a  complex  $g$-dimensional  vector  space, so the moduli space
$\cH$ of pairs (Riemann surface $C$, holomorphic $1$-form $\omega$ on
$C$)  is  a  total  space  of  complex  $g$-dimensional \textit{Hodge
bundle}  over  the moduli space of complex curves $\cM_g$. The moduli
space   $\cH$   is   called   \textit{the  moduli  space  of  Abelian
differentials}.  A holomorphic $1$-form $\omega$ on a Riemann surface
of  genus $g$ has $2g-2$ zeroes counting multiplicities. For example,
a  holomorphic  1-form  on  a Riemann surface of genus two might have
either  two  simple  zeroes or a single zero of order $2$. The moduli
space    $\cH$    gets    naturally    stratified   by   the   strata
$\cH(m_1,\dots,m_n)$  of  pairs  (Riemann  surface  $C$,  holomorphic
$1$-form  $\omega$  on  $C$  with  zeroes of orders $m_1,\dots,m_n$),
where  $m_1+\dots+m_n=2g-2$. For example, in genus two there are only
two  such  strata,  namely, the \textit{principal} stratum $\cH(1,1)$
and the stratum $\cH(2)$.

Note  that,  in  general,   individual strata do not carry any bundle
structure over $\cM_g$. The dimension of the stratum is computed as
$$
\dim_{\C}\cH(m_1,\dots,m_n)=2g+n-1\,.
$$
In  particular,  the  dimension of the minimal stratum $\cH(2g-2)$ is
$2g$  while  the  dimension  of the moduli space of curves $\cM_g$ is
$3g-3$.

The strata are not necessarily connected: they might have up to three
connected  components.  In  genera  greater than or equal to four the
classifying  invariants are the \textit{parity of the spin structure}
applicable  to  strata  $\cH(2k_1,\dots,2k_n)$  for which all entries
$m_i=2k_i$  are even, and \textit{hyperellipticity} applicable to two
special   strata,  namely,  to  $\cH(g-1,g-1)$  and  to  $\cH(2g-2)$,
see~\cite{Kontsevich:Zorich}.

Denote by $P_1,\dots,P_n\in C$ the points of the Riemann surface $C$,
where   the   holomorphic   1-form   has  zeroes.  The  vector  space
$H^1(C,\{P_1,\dots,P_n\};\C)$ provides canonical local \textit{period
coordinates}   in   the  stratum  $\cH(m_1,\dots,m_n)$.  The  lattice
$H^1(C,\{P_1,\dots,P_n\};\Z\oplus  i\Z)$  defines  a canonical volume
element  $d\nu$ in the stratum. This volume element induces a natural
volume  element on the real hypersurface $\cH_1(m_1,\dots,m_n)\subset
\cH(m_1,\dots,m_n)$ given by the equation
$$
\frac{i}{2}\int_C \omega\wedge\bar\omega =1\,.
$$

The   group   $\GL$  acts  naturally  on  every  stratum:  in  period
coordinates
$$
H^1(C,\{P_1,\dots,P_n\};\C)\simeq
H^1(C,\{P_1,\dots,P_n\};\R)\otimes\R^2
$$
$\GL$  acts  on  the  second term in the tensor product. The subgroup
$\SL$  preserves  the  volume  element  $d\nu_1$.  By the fundamental
theorem  of H.~Masur~\cite{Masur} and W.~Veech~\cite{Veech} the total
volume  $\nu_1(\cH_1(m_1,\dots,m_n))$  of  any stratum is finite, and
the  action  of  $\SL$ is ergodic on every connected component of the
stratum with respect to the invariant measure $d\nu_1$.

A  holomorphic  1-form  $\omega$  defines  a flat metric on a Riemann
surface  $C$  with  conical  singularities at the zeroes of $\omega$.
Such flat metrics appear, for example, by unfolding rational billiard
tables.  It  is  very useful to apply the following technology in the
study  of  behavior  of  geodesics  on  flat  surfaces  of this kind.
Consider  the  $\GL$-orbit  of  the point $(C,\omega)$ in the ambient
stratum   of   Abelian  differentials  and  find  the  orbit  closure
$\overline{\GL\cdot(C,\omega)}$.  The  geometry  of the orbit closure
provides you  with  important information on geometry and dynamics of
the  original  flat  surface. The theorem of H.~Masur and of W.~Veech
tells  that  for almost all flat surfaces in each stratum, such orbit
closure  is the entire ambient connected component of the stratum. In
particular,  several  important  characteristics  of  individual flat
surfaces    (such    as    Siegel--Veech   constants   discussed   in
section~\ref{s:Siegel:Veech})  are  expressed in terms of the volumes
of the connected components of the strata.

For  a  deeper  presentation  of these ideas we recommend a short and
extremely accessible introductory article~\cite{Wright:BAMS} and more
detailed            survey           papers~\cite{Masur:Tabachnikov},
\cite{Wright:survey}, \cite{Zorich:Houches}.

\section{Conjectural           asymptotics           for          the
volumes of the strata of Abelian   differentials}

Let  $m=(m_1,\dots,m_n)$  be an unordered partition of an even number
$2g-2$,  i.e.,  let  $|m|=m_1+\dots+m_n=2g-2$. Denote by $\Pi_{2g-2}$
the set of all such partitions.

\begin{MainConjecture}
\label{mc:vol}
For any $m\in\Pi_{2g-2}$ one has
\begin{equation}
\label{eq:asymptotic:formula:for:the:volume}
\Vol\cH_1(m_1,\dots,m_n)=\cfrac{4}{(m_1+1)\cdot\dots\cdot(m_n+1)}
\cdot (1+\varepsilon_1(m)),
\end{equation}
where
$$
\lim_{g\to\infty} \max_{m\in\Pi_{2g-2}} |\varepsilon_1(m)| = 0.
$$
\end{MainConjecture}

\begin{Remark}
It was  proved  in~\cite{Eskin:Okounkov}  that for any partition
$m\in\Pi_{2g}$ one has
$$
\cfrac{\Vol\cH_1(m_1,\dots,m_n)}{\pi^{2g}}\in\Q.
$$
Note,  however, that  the  number $\pi$ is  not  involved in  the
asymptotic   formula~\eqref{eq:asymptotic:formula:for:the:volume}
anymore.
\end{Remark}

\begin{Remark}
We  use  the normalization of volumes from~\cite{Eskin:Masur:Zorich}.
In   particular,   all   the   zeroes  are  ``named''.  In  notations
of~\cite{Eskin:Okounkov} the formula for the volume should be read as
\begin{equation*}
\label{eq:vol:equals:c}
\Vol(\cH_1(m))=2\cdot{\bf c}(m+\vec{1})
\end{equation*}
\end{Remark}

We  conjecture  that  convergence~\eqref{eq:asymptotic:formula:for:the:volume}
is sufficiently fast.

\begin{Conjecture}
\label{conj:strong}
There exists a universal constant $C$ such that for all $g$
and all $m\in\Pi_g$ one has
\begin{equation}
\label{eq:strong:epsilon1}
|\epsilon_1(m)|\le \frac{C}{\sqrt{g}}\,.
\end{equation}
\end{Conjecture}

Consider any fixed
collection    of    positive    integers    $m_1,..,m_n$,   and   let
$|m|=m_1+\dots+m_n$.  For any integer $g$ such that $2g-2>|m|$ define
$m(g)$ as $m$ completed by $2g-2-|m|$ entries $1$,
\begin{equation}
\label{eq:add:ones}
m(g):=(m_1,\dots,m_n,\underbrace{1,\dots,1}_{2g-2-|m|})
\end{equation}
We expect that the convergence is faster
than~\eqref{eq:strong:epsilon1} for strata~\eqref{eq:add:ones} when
$g\to+\infty$.

\begin{Remark}
In  the  conjectures  above  we  considered
the asymptotics with respect to the genus $g$ of the stratum.
We  could also consider the asymptotics with respect to
the dimension of the stratum. Since the complex dimensions
of  the  strata  in a given  genus  $g$  vary  from $2g$ to $4g-3$ the
corresponding asymptotics are, basically, equivalent.
\end{Remark}

Note  that  some  strata  of  Abelian  differentials  contain several
connected  components~\cite{Kontsevich:Zorich}. The conjecture on the
asymptotics  of  the  volumes  is stated for the total volume of such
strata.  However,  it  can  be  easily  specified  for the individual
components.  The  result  in~\cite{Athreya:Eskin:Zorich}  provides the
exact  value  for  the  hyperelliptic connected components (and, more
generally, for all hyperelliptic loci), namely:

\begin{align*}
\Vol\cH^{hyp}_1(2g-2)&=\cfrac{2\pi^{2g}}{(2g+1)!}\cdot
\cfrac{(2g-3)!!}{(2g-2)!!} \sim
\cfrac{1}{\pi^2 g}\left(\frac{\pi e}{2g+1}\right)^{2g+1}\,.
\\
\Vol\cH^{hyp}_1(g-1,g-1)&=\cfrac{4\pi^{2g}}{(2g+2)!}\cdot
\cfrac{(2g-2)!!}{(2g-1)!!} \sim
\cfrac{1}{\pi^2 g}\left(\frac{\pi e}{2g+2}\right)^{2g+2}\,.
\end{align*}

The  above  formulae  show  that  the  volume  of  the  hyperelliptic
connected components is negligible in comparison with the conjectural
volume~\eqref{eq:asymptotic:formula:for:the:volume}  of  the  stratum
when genus is sufficiently large.

\begin{Conjecture}
For any $k\in\Pi_{g-1}$ one has
\begin{equation}
\label{eq:even:over:odd}
\cfrac{\Vol\cH^{even}_1(2k_1,\dots,2k_n)}
{\Vol\cH^{odd}_1(2k_1,\dots,2k_n)}
= (1+\varepsilon_2(k)),
\end{equation}
where
$$
\lim_{g\to\infty} \max_{k\in\Pi_{g-1}} \varepsilon_2(k) = 0.
$$
\end{Conjecture}

This  Conjecture  might  be  approached  by a crude estimate since we
expect     that     the     difference     between
$\Vol\cH^{even}_1(2k_1,\dots,2k_n)$                               and
$\Vol\cH^{odd}_1(2k_1,\dots,2k_n)$ is much smaller than their values.

\begin{figure}[htb]
%
   %
\includegraphics{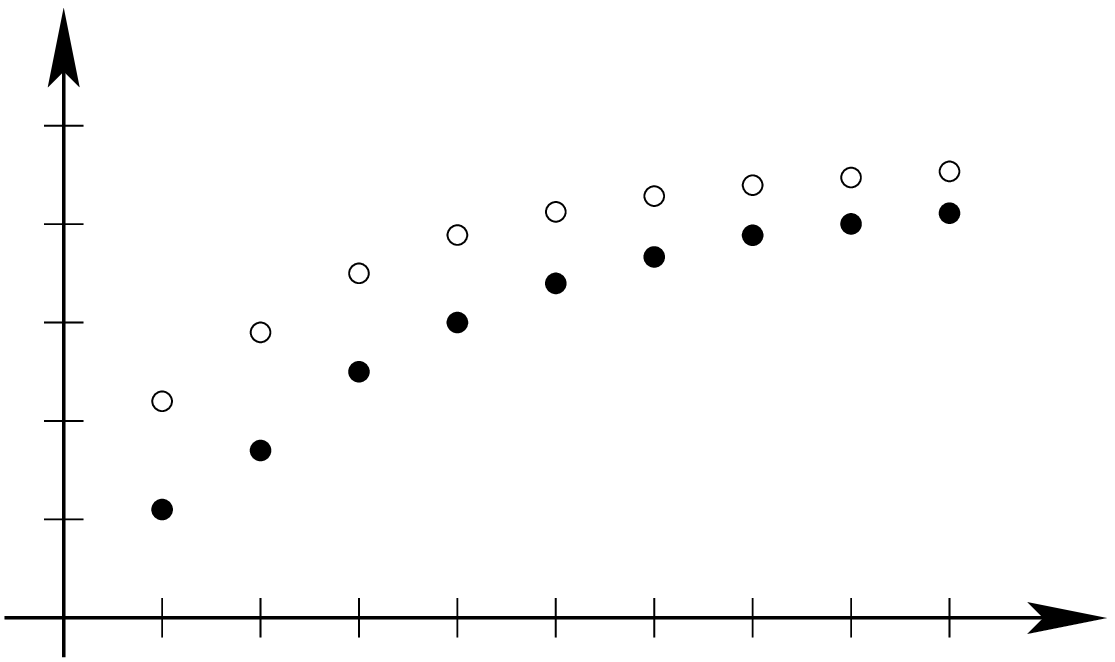}
\includegraphics{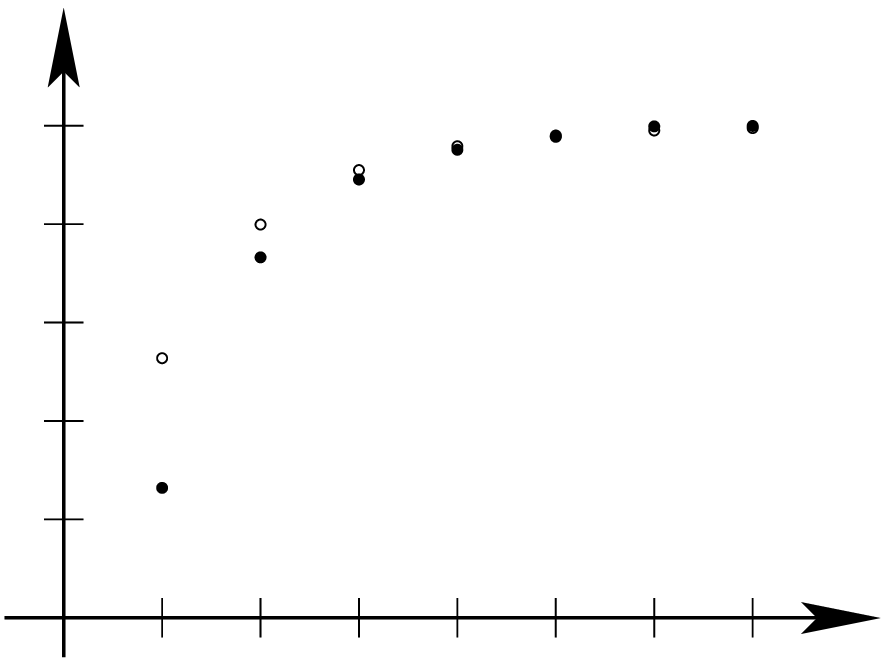}
\begin{picture}(0,0)(0,0)
\begin{picture}(0,0)(0,0)
\put(-142,-97){$2$}
\put(-128,-97){$3$}
\put(-114,-97){$4$}
\put(-100,-97){$5$}
\put(-86,-97){$6$}
\put(-71,-97){$7$}
\put(-57,-97){$8$}
\put(-43,-97){$9$}
\put(-32,-97){$10$}
\put(-20,-80){$g$}
\put(-172,-73){\small $0.6$}
\put(-172,-59){\small $0.7$}
\put(-172,-44){\small $0.8$}
\put(-172,-30){\small $0.9$}
\put(-172,-16){\small $1.0$}
\put(-145,-5){$1+\epsilon_1(m)$}
\end{picture}
\begin{picture}(0,0)(-186,0)
\put(-142,-97){$4$}
\put(-128,-97){$5$}
\put(-114,-97){$6$}
\put(-100,-97){$7$}
\put(-86,-97){$8$}
\put(-71,-97){$9$}
\put(-57,-97){$10$}
\put(-50,-80){$g$}
\put(-172,-73){\small $0.6$}
\put(-172,-59){\small $0.7$}
\put(-172,-44){\small $0.8$}
\put(-172,-30){\small $0.9$}
\put(-172,-16){\small $1.0$}
\put(-145,-5){$1+\epsilon_2(2k)$}
\end{picture}
\end{picture}
\vspace{100bp}
\caption{
\label{fig:epsilon1and2}
Maximum and minimum of $(1+\epsilon_1(m))$ over all
strata in fixed genus $g=2,\dots, 10$ on the left.
Maximum and minimum of $(1+\epsilon_2(2k))$ over all
admissible strata in fixed genus $g=4,\dots, 10$ on the
right.
   }
\end{figure}

\noindent  \textbf{Numerical  evidence.}  Actually, in the context of
the  Main  Conjecture~\ref{mc:vol},  the infinity does not seem to be
very far. We present here the graphs of exact numerical data for
$$
\max_{m\in\Pi_{2g-2}} \left(1+\epsilon_i(m)\right) \text{ and }
\min_{m\in\Pi_{2g-2}} \left(1+\epsilon_i(m)\right),
\ i=1,2\,.
$$
for genera $g\le 10$. The reader can see that at least for this range
of genera, one always has $\epsilon_i(m)<0$ and that both maximal and
minimal  value  of  $\epsilon_i(m)$  over all strata in a given genus
monotonously  tend to zero as genus grows.

For any given genus in the
range  $2\le  g\le  10$  the  minimum  of $|\epsilon_1(m)|$ is always
attained  on the principal stratum $\cH(1^{2g-2})$ and the maximum of
$|\epsilon_1(m)|$ is always attained on the (non connected for $g>2$)
stratum  $\cH(2g-2)$.  In  other words (at least for this range of
genera) the convergence in~\eqref{eq:asymptotic:formula:for:the:volume}
is the best for the principal stratum and the
worst for the disconnected stratum $\cH(2g-2)$.

We have $\epsilon_2(2k)<0$ for all partitions $k$
in genera in the range from $4$ to $10$.
The   minimum   of   $|\epsilon_2(2k)|$ (best convergence in~\eqref{eq:even:over:odd})
over  all  partitions  $k\in
\Pi_{g-1}$  for any fixed genus in this range is always
attained   on   the   stratum   $\cH(2^{g-1})$  and  the  maximum  of
$|\epsilon_2(2k)|$ (worst convergence  in~\eqref{eq:even:over:odd})
is always attained on the stratum $\cH(2g-2)$.

\medskip

\noindent
\textbf{State of the art.}
The  Main  Conjecture~\ref{mc:vol}  was
recently  proved by D.~Chen, M.~M\"oller, and
D.~Zagier~\cite{ChenMollerZagier}
for  the  principal  stratum $\cH(1^{2g-2})$.    They use   and   develop
quasimodularity  properties  of a new
\textit{all genera generating function}
for volumes.

\section{Conjectural   universality  of  the  Siegel--Veech  constant
for strata of Abelian   differentials of large genus}
\label{s:Siegel:Veech}

\begin{MainConjecture}
\label{mc:SV}
For  all  nonhyperelliptic  connected  components $\cH$ of all strata
$\cH(m)$ of Abelian differentials, where $m\in\Pi_{2g-2}$, one has
$$
\lim_{g\to\infty} \carea(\cH) = \frac{1}{2}\,.
$$
\end{MainConjecture}

\begin{ConditionalCorollary}
Applying  the formula from~\cite{Eskin:Kontsevich:Zorich} for the sum
of  the Lyapunov exponents of the Hodge bundle over the Teichm\"uller
geodesic  flow  in  any  nonhyperelliptic  connected  component  of a
stratum $\cH_1(m)$ of Abelian differentials we get
\begin{equation}
\label{eq:general:sum:of:exponents:for:Abelian}
\lambda_1 + \dots + \lambda_g
\ = \
\cfrac{1}{12}\cdot\sum_{m_i\in m} \cfrac{m_i(m_i+2)}{m_i+1}
\ +\ \frac{\pi^2}{6}+\epsilon_3(m)\,,
\end{equation}
where $\lim_{g\to+\infty}\max_{m\in\Pi_{2g-2}}|\epsilon_3(m)|\to 0$.
\end{ConditionalCorollary}

Recall that the original definition of the Siegel--Veech constant
comes from counting of closed geodesics. Regular simple closed
geodesics on translation surfaces appear in families; every such
family fills a maximal flat cylinder having at least one conical
singularity on each of the two boundary components. Consider a
translation surface of unit area, and let $N_{\mathit{area}}(S,L)$ be
the weighted number of regular simple closed geodesics of length at
most $L$, where the weight is the area of the associated maximal flat
cylinder. Morally, we pretend that the geodesic is thick, and the
weight measures its thickness. By the result of A.~Eskin and
H.~Masur~\cite{Eskin:Masur}, for almost all flat surfaces $S$ of unit
area in any connected component $\cH$ of any stratum of Abelian
differentials, the weighted number $N_{\mathit{area}}(S,L)$ of
regular simple closed geodesics has quadratic asymptotics
$$
N_{\mathit{area}}(S,L)=\carea(\cH)\cdot\pi L^2 +o(L^2)\ \text{ as }L\to+\infty\,.
$$
Morally, Conjecture~\ref{mc:SV} claims that for a random translation
surface of sufficiently large genus, the asymptotics
of $N_{\mathit{area}}(S,L)$ is, basically, universal
no matter what is the genus and what are the number and
types of conical singularities of $S$.

\begin{Remark}
Due to a misunderstanding in the normalization, the conjectural value
of $\carea$ indicated in~\cite{Chen} differs from the actual one by a
constant factor.
\end{Remark}

\begin{Remark}
The  hyperelliptic  connected  components  excluded in the statements
above        exhibit        completely       different       behavior
(see~\cite{Eskin:Kontsevich:Zorich}):
\begin{align*}
\carea(\cH^{hyp}_1(2g-2))&=\cfrac{(2g+1)g}{\pi^2(2g-1)}
\\
\carea(\cH^{hyp}_1(g-1,g-1))&=\cfrac{(2g+1)(g+1)}{\pi^2(2g)}\,.
\end{align*}
\end{Remark}
\medskip

\begin{figure}[htb]
%
   %
\includegraphics{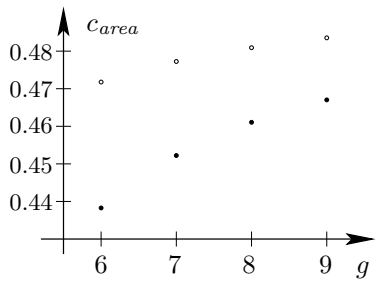}
\begin{picture}(0,0)(-95,0)
\put(-142,-97){$6$}
\put(-114,-97){$7$}
\put(-86,-97){$8$}
\put(-57,-97){$9$}
\put(-43,-97){$g$}
\put(-174,-73){\small $0.44$}
\put(-174,-59){\small $0.45$}
\put(-174,-44){\small $0.46$}
\put(-174,-30){\small $0.47$}
\put(-174,-16){\small $0.48$}
\put(-145,-5){$\carea$}
\end{picture}
\vspace{100bp}
\caption{
\label{fig:minmax:carea}
Maximum and minimum of $\carea$ over all
non hyperelliptic connected components of all
strata in fixed genus $g=6,\dots, 9$.
   }
\end{figure}

\noindent  \textbf{Numerical  evidence.}
There  are  two  sources  of strong numerical evidence supporting the
above   conjecture.   On  the  one  hand  we  have  exact  values  of
$\carea(\cH)$  for all connected components of all strata up to genus
$g=9$. For this range of genera we have $\carea(\cH)<\frac{1}{2}$ for
all non hyperelliptic components of all strata.
The maximal value of $\carea$ (closest to $\frac{1}{2}$)
in this range of genera is
attained  for  the  principal  stratum; the minimal value --- for the
component $\cH^{odd}(2g-2)$, and both are already close to $\frac{1}{2}$,
see~Figure~\ref{fig:minmax:carea}.

Numerical  evidence  which  covers strata of larger genera comes from
numerical  simulations  for  the  sums of the Lyapunov exponents. For
those   strata   where  we  have  exact  values  of  $\carea$,  these
simulations  give  up to five digits of precision, which gives a hope
that they are sufficiently reliable in higher genera. Also,
from~\cite{Eskin:Kontsevich:Zorich} we know the exact value of
$\carea(\cH^{hyp})$ for all genera, and it matches the numerical
simulations based on evaluation of Lyapunov exponents. These kind
of numerical simulations do not seem to be sensitive to a particular
stratum.
\medskip

\noindent  \textbf{State of the art.}
The  Main  Conjecture~\ref{mc:SV}  was  recently  proved  by D.~Chen,
M.~M\"oller,  and D.~Zagier~\cite{ChenMollerZagier} for the principal
stratum  $\cH(1,\dots,1)$. The proof is also based on quasimodularity
properties  of a new generating function which the authors
elaborate for $\carea$. The method is applicable to
strata~\eqref{eq:add:ones} as $g\to+\infty$.

We have an alternative \textit{conditional} straightforward proof  of
the  Main Conjecture~\ref{mc:SV} for the principal stratum
$\cH(1,\dots,1)$ relying on its strong
version~\eqref{eq:strong:epsilon1} proved  in~\cite{ChenMollerZagier}
and on the original formula from~\cite{Eskin:Masur:Zorich} for
$\carea$ in terms of  the polynomial in volumes. The same method
works for  strata~\eqref{eq:add:ones} as $g\to+\infty$ as soon as the
Main Conjecture     is     proved    in    the    strong
form~\eqref{eq:strong:epsilon1}  for  the initial stratum $\cH(m(g))$
and for all analogous strata obtained  by  replacing  any
subcollection  of $m_i$ by all possible partitions of $m_i-1$.

\subsection*{Acknowledgements}
We are grateful MPIM in Bonn for constant hospitality along many
years of development of this project. We thank D.~Chen, M.~M\"oller
and D.~Zagier for their interest to the conjectures and for
important conversations.


\end{document}